\documentclass[fleqn,twoside]{article}

\usepackage{amsmath}
\usepackage{amsfonts}
\usepackage{amsthm}
\usepackage{amscd}
\usepackage{espcrc2}

\usepackage{graphicx}
\usepackage[figuresright]{rotating}

\newtheorem{thm}{Theorem}[section]
\newcommand{\B}{\mathcal B}
\newcommand{\C}{\mathcal C}
\newcommand{\E}{\mathcal E}
\newcommand{\ep}{\epsilon}
\newcommand{\zt}{\zeta}
\newcommand{\si}{\sigma}
\newcommand{\De}{\Delta}
\newcommand{\Si}{\Sigma}
\newcommand{\Li}{\operatorname{Li}}
\newcommand{\id}{\operatorname{id}}
\newcommand{\card}{\operatorname{card}}
\newcommand{\Sym}{\operatorname{Sym}}
\newcommand{\QSym}{\operatorname{QSym}}


\title{The Hopf algebra structure of multiple harmonic sums}

\author{Michael E. Hoffman\address{Dept. of Mathematics, U. S. Naval Academy, 
        Annapolis, MD 21402}}
       
\begin{document}

\begin{abstract}
Multiple harmonic sums appear in the perturbative computation of 
various quantities of interest in quantum field theory.  
In this article we introduce
a class of Hopf algebras that describe the structure of such sums,
and develop some of their properties that can be exploited in
calculations.
\vspace{1pc}
\end{abstract}

\maketitle

\section{MULTIPLE HARMONIC SUMS}

As discussed in the introduction of \cite{B03}, multiple harmonic sums 
occur in perturbative higher-order calculations of quantum field theory.
Let $I=(i_1,i_2,\dots,i_k)$ be a sequence of positive integers.
For positive integers $n$, we define the multiple harmonic sums
\begin{multline}
A_I(n;x_1,x_2\dots,x_k)=\\
\sum_{n\ge n_1>n_2>\cdots>n_k\ge 1}
\frac{x_1^{n_1}x_2^{n_2}\cdots x_k^{n_k}}{n_1^{i_1}n_2^{i_2}\cdots n_k^{i_k}} 
\label{Adef}
\end{multline}
and
\begin{multline}
S_I(n;x_1,x_2\dots,x_k)=\\
\sum_{n\ge n_1\ge n_2\ge\cdots\ge n_k\ge 1}
\frac{x_1^{n_1}x_2^{n_2}\cdots x_k^{n_k}}{n_1^{i_1}n_2^{i_2}\cdots n_k^{i_k}} 
\label{Sdef}
\end{multline}
associated with $I$ (Note that the only difference between
(\ref{Adef}) and (\ref{Sdef}) is in the inequalities in the summation 
variables).
Both types of sums appear in \cite{MUW02} and \cite{W03}, with
a slightly different notation ($Z$ is used in place of our $A$).
\par
If $i_1>1$ and $x_1=x_2=\dots=x_k=1$, the sums (\ref{Adef}) and (\ref{Sdef}) 
converge as $n\to\infty$, giving the well-known multiple zeta values 
\cite{H92,Z94,H97,Z03}:
$$
\zt(i_1,\dots,i_k)=A_{(i_1,\dots,i_k)}(\infty;1,\dots,1),
$$
which also occur in some perturbative QFT calculations \cite{BK97}.  
More generally,
(\ref{Adef}) and (\ref{Sdef}) converge as $n\to\infty$ when $|x_i|=1$ 
for all $i$ and $i_1x_1\ne 1$.  
The quantities 
$$
\Li_I(x_1,\dots,x_n)=A_I(\infty;x_1,\dots,x_k)
$$
are called 
multiple polylogarithms \cite{BB01}:  they generalize the classical 
polylogarithm $\Li_n(x_1)=A_{(n)}(\infty;x_1)$.
\par
It is immediate from the defining equations (\ref{Adef}) and (\ref{Sdef}) 
that the sums $S_I$ can be written in terms of the $A_I$.  
To state the 
relation precisely, let $\C(n)$ be the set of compositions of $n$,
i.e., ordered sequences $(i_1,\dots,i_k)$ of positive integers
with $i_1+\dots+i_k=n$.  If $I=(i_1,\dots,i_k)$ is a composition of 
$n$ and $J=(j_1,\dots,j_p)$ is a composition of $k$, then there
is a composition $J\circ I$ of $n$ given by
\begin{multline*}
(i_1+\dots+i_{j_1},i_{j_1+1}+\dots+i_{j_1+j_2},\dots,\\
i_{k-j_p+1}+\dots+i_k)
\end{multline*}
(cf. \cite[p. 52]{H00}).  Also, compositions act on argument strings:
given $J=(j_1,\dots,j_p)\in\C(k)$ and a string $X=(x_1,\dots,x_k)$ 
of length $k$, we have
\begin{multline*}
J(X)=\\
(x_1\cdots x_{j_1},x_{j_1+1}\cdots x_{j_1+j_2},\dots,
x_{k-j_p+1}\cdots x_k) .
\end{multline*}
Then the relation between sums of types (\ref{Adef}) and (\ref{Sdef}) 
is given by
\begin{equation}
S_I(n;X)=\sum_{J\in\C(k)}A_{J\circ I}(n;J(X))
\label{StoA}
\end{equation}
for any $I=(i_1,\dots,i_k)$ and $X=(x_1,\dots,x_k)$.
M\"obius inversion can be applied to (\ref{StoA}) to obtain
\begin{equation}
A_I(n;X)=\sum_{J\in\C(k)}(-1)^{\ell(J)-k}S_{J\circ I}(n;J(X)),
\label{AtoS}
\end{equation}
where $\ell(J)$ is the number of parts of $J$.  In fact, there is
a deeper relation between the $S_I$ and the $A_I$ than
the essentially trivial equations (\ref{StoA}) and (\ref{AtoS}):
the two types of sums are dual in the sense that they have,
up to signs, the same algebraic properties.  
In the case where the arguments $x_i$ are roots of unity, we
can formalize the algebra of such sums using the family
of Hopf algebras described in the next section.  
(Making the arguments roots of unity seems to capture many cases 
of physical interest:  see, e.g., \cite{B99}.)

\section{THE EULER ALGEBRA}

We recall from \cite{H00} the construction of the Euler algebra of 
index $r$, where $r$ is a positive integer.  
We start with noncommuting symbols (or ``letters'') $z_{i,j}$,
where $i,j$ are integers with $i$ positive and $0\le j\le r-1$.
Let $\E_r$ be the complex vector space generated by words in the
$z_{i,j}$ (including the empty word, denoted by 1).  For such a
word $w=z_{i_1,j_1}z_{i_2,j_2}\cdots z_{i_k,j_k}$, we define
the degree of $w$ to be $|w|=i_1+\dots+i_k$ (and call $\ell(w)=k$ 
the length of $w$.)
Now we define a multiplication $*$ on $\E_r$
as follows.
We require $1*w=w*1=w$ for all words $w$, and 
\begin{multline}
z_{i,j}w_1*z_{m,n}w_2=z_{i,j}(w_1*z_{m,n}w_2)+\\
z_{m,n}(z_{i,j}w_1*w_2)
+z_{i+m,j+n}(w_1*w_2)
\label{mult1}
\end{multline}
for any words $w_1,w_2$:  here the addition in the second subscript
is to be understood mod $r$.  For example, when $r=3$ 
\begin{multline*}
z_{1,1}*z_{1,2}z_{2,1}=z_{1,1}z_{1,2}z_{2,1}+z_{1,2}(z_{1,1}*z_{2,1})
+z_{2,0}z_{2,1}\\
=z_{1,1}z_{1,2}z_{2,1}+z_{1,2}z_{1,1}z_{2,1}+z_{1,2}z_{2,1}z_{1,1}+\\
z_{1,2}z_{3,2}+z_{2,0}z_{2,1} .
\end{multline*}
Since each of the parenthesized products on the right-hand side 
of (\ref{mult1}) has total length less than the left-hand side, 
equation (\ref{mult1}) gives an inductive definition of $*$ on $\E_r$.  
As shown in \cite{H00}, $(\E_r,*)$ is a commutative, associative
graded algebra over $\mathbb C$.  In fact, $(\E_r,*)$ is a polynomial
algebra.  To describe the generators, we first assume that the
letters $z_{i,j}$ are totally ordered, and extend this order lexicographically
to words.  A word $w$ is called Lyndon if it is smaller than any of its
proper right factors ($v\ne 1$ is a proper right factor of $w$ if
$w=uv$ for $u\ne 1$).  Then we have the following result
\cite[Theorem 2.6]{H00}.
\begin{thm} For positive integer $r$, $(\E_r,*)$ is the polynomial
algebra on the Lyndon words.
\end{thm}
\par\noindent
\emph{Remark.} From the discussion in \cite[Example 2]{H00}, the
number of Lyndon words of degree $n$ in $\E_r$ is
$$
\frac1{n}\sum_{d|n}\mu\left(\frac{n}{d}\right)(r+1)^d ,
$$
where the sum is over divisors of $n$ and $\mu$ is the M\"obius
function on the integers.
\par
In the case $r=1$, $(\E_r,*)$ is the algebra $\QSym$ of quasi-symmetric
functions as defined by Gessel \cite{G84}.  For a description of 
$\QSym$ and its relation to multiple harmonic sums with unit arguments,
see \cite{H04}.  Note that $\QSym$ contains the well-known algebra $\Sym$ 
of symmetric functions:  in fact, $\Sym$ can
be imbedded in any $\E_r$ by sending the elementary symmetric function
$e_i$ to $z_{1,0}^i$.
\par
We can also define a coalgebra structure on $\E_r$ as follows.
The counit $\ep:\E_r\to\mathbb C$ is given by
$$
\ep(1)=1,\quad \ep(w)=0\quad\text{for $|w|>0$}
$$
and the coproduct $\De:\E_r\to\E_r\otimes\E_r$ by
\begin{multline*}
\De(z_{i_1,j_1}z_{i_2,j_2}\cdots z_{i_k,j_k})=\\
\sum_{p=0}^k z_{i_1,j_1}\cdots z_{i_p,j_p}\otimes z_{i_{p+1},j_{p+1}}\cdots
z_{i_k,j_k} .
\end{multline*}
Then Theorem 3.1 of \cite{H00} says that $(\E_r,*,\De)$ is a graded
connected Hopf algebra.  
\par
Since the $*$-product is commutative, the antipode $S$ of $(\E_r,*,\De)$ 
is an involution, i.e., an algebra automorphism with $S^2=\id$ (see,
e.g., \cite[Theorem III.3.4]{K95}).
As shown in \cite[Theorem 3.2]{H00}, 
there are two (not obviously identical) formulas for $S$.  
Our first formula for $S$ involves iterated products:  
\begin{equation}
S(w)=\sum_{w_1w_2\cdots w_k=w}(-1)^k w_1*w_2*\cdots*w_k,
\label{ant1}
\end{equation}
where the sum is over all decompositions of $w$ into (nonempty) 
subwords $w_1,\dots,w_k$.
\par
For the second formula, we shall introduce some more notation.
Given a string $a_1,\dots, a_n$ of letters, let $[a_1,\dots,a_n]$
be the letter obtained by adding all the subscripts (where of 
course the addition in the second subscript is mod $r$).
Then the multiplication rule (\ref{mult1}) can be written
\begin{multline}
aw*bv=\\
a(w*bv)+b(aw*v)+[a,b](w*v)
\label{mult2}
\end{multline}
for letters $a,b$ and words $w,v$.  
As above, let $\C(n)$ be the set of compositions of $n$.
Then $(i_1,\dots,i_k)\in\C(n)$ acts on a word $w=a_1\cdots a_n$ of 
length $n$ as follows:
\begin{multline*}
(i_1,\dots,i_k)[w]=\\
[a_1,\dots,a_{i_1}][a_{i_1+1},\dots,a_{i_1+i_2}]\cdots
[a_{n-i_k+1},\dots,a_n].
\end{multline*}
Our second formula for the antipode can be written in terms of this
action:
\begin{equation}
S(w)=(-1)^n\sum_{I\in\C(n)}I[a_na_{n-1}\cdots a_1] 
\label{ant2}
\end{equation}
for words $w=a_1\cdots a_n$ of length $n$.
\par
Now let $R:\E_r\to\E_r$ be the linear function that reverses words,
i.e.,
$$
R(a_1a_2\cdots a_n)=a_n\cdots a_2a_1.
$$
The following result can be proved by induction on word length
(see \cite[Theorem 9]{Z03}).
\begin{thm}
The function $R:\E_r\to\E_r$ is a $*$-automorphism.
\end{thm}
Since clearly $\De\circ R=(R\otimes R)\circ\De$, $R$ is evidently
an automorphism of the Hopf algebra $(\E_r,*,\De)$.
\par
The action of compositions on words of $\E_r$ can be used to put
a partial order $\preceq$ on words as follows.  For a word 
$w=a_1\cdots a_n$ of length $n$, set $v\preceq w$ if $v=I[w]$ for
some $I\in\C(n)$ (Note that in this case $\ell(v)\le \ell(w)$ and
$|v|=|w|$).  Define
\begin{equation}
\overline w=\sum_{v\preceq w} v =\sum_{J\in\C(\ell(w))}J[w]
\label{bardef}
\end{equation}
for words $w$ of $\E_r$.  
\par
Our second formula for the antipode can now be written 
$$
RS(w)=SR(w)=(-1)^{\ell(w)}\overline w
$$
for any word $w$ of $\E_r$.  Equating the two formulas for the
antipode, we have
\begin{multline*}
\sum_{w_1w_2\cdots w_k=w}(-1)^kw_1*w_2*\cdots*w_k=\\
(-1)^{\ell(w)}\sum_{v\preceq w}R(v),
\end{multline*}
or, after applying $R$ to both sides,
\begin{equation}
\overline w=\sum_{w_1\cdots w_k=R(w)}(-1)^{\ell(w)-k}w_1*\cdots*w_k .
\label{barw2w}
\end{equation}
Now apply $RS$ to both sides of equation (\ref{barw2w}):
\begin{equation}
w=\sum_{w_1\cdots w_k=R(w)}(-1)^{\ell(w)-k}\overline w_1*\cdots*
\overline w_k .
\label{w2barw}
\end{equation}
Also, since $\overline w=(-1)^{\ell(w)}SR(w)$ and $SR$ is an automorphism
of the Hopf algebra $(\E_r,*,\De)$, we can work with the vector space
basis $\overline w$ just as well as with the basis consisting of the 
words $w$:  the only difference is that the inductive rule (\ref{mult2}) 
for the $*$-product is replaced by 
$$
\overline{aw}*\overline{bv}=\overline{a(w*bv)}
+\overline{b(aw*v)}-\overline{[a,b](w*v)} .
$$

\section{RELATION TO MULTIPLE HARMONIC SUMS}

Now we relate the Hopf algebras $\E_r$ to the multiple harmonic
sums.
For fixed $n$, define a linear map $\rho_n:\E_r\to\mathbb C$ by
$$
\rho_n(z_{i_1,j_1}\cdots z_{i_k,j_k})=
A_{(i_1,\dots,i_k)}(n;\ep^{j_1},\dots,\ep^{j_k})
$$
where $\ep=e^{\frac{2\pi i}{r}}$.  We have the following result.
\begin{thm}
The function $\rho_n$ is a homomorphism of $(\E_r,*)$ into $\mathbb C$.
\end{thm}
\begin{proof} $\rho_n$ is the composition of the homomorphism
$\phi_n:\E_r\to\mathbb C[t_1,\dots,t_n]$ of \cite[Theorem 7.1]{H00}
with the homomorphism $\mathbb C[t_1\dots,t_n]\to\mathbb C$ sending
$t_i$ to $1/i$, $1\le i\le n$.
\end{proof}
Comparing the definition (\ref{bardef}) of $\overline w$
with equation (\ref{StoA}), it is evident that
$$
\rho_n(\overline{z_{i_1,j_1}\cdots z_{i_k,j_k}})=
S_{(i_1,\dots,i_k)}(n;\ep^{j_1},\dots,\ep^{j_k}) .
$$
Henceforth we shall write $A_w(n)$ for $\rho_n(w)$ and $S_w(n)$
for $\rho_n(\overline w)$ for words $w$ of $\E_r$, using the word
to code for both the exponents and roots-of-unity arguments.
In this notation, equation (\ref{StoA}) is
\begin{equation}
S_w(n)=\sum_{u\preceq w}A_u(n) .
\label{stoa}
\end{equation}
Now applying $\rho_n$ to equations (\ref{barw2w}) and 
(\ref{w2barw}) gives respectively
\begin{multline}
S_w(n)=\\
\sum_{w_1\cdots w_k=R(w)}(-1)^{\ell(w)-k}A_{w_1}(n)\cdots
A_{w_k}(n)
\label{s2proda}
\end{multline}
and
\begin{multline}
A_w(n)=\\
\sum_{w_1\cdots w_k=R(w)}(-1)^{\ell(w)-k}S_{w_1}(n)\cdots
S_{w_k}(n) .
\label{a2prods}
\end{multline}
By equating the right-hand sides of equations (\ref{stoa})
and (\ref{s2proda}), one obtains
\begin{multline*}
A_w(n)+(-1)^{\ell(w)}A_{R(w)}(n)=\\
\sum_{\substack{w_1\cdots w_k=R(w)\\k>1}}
(-1)^{\ell(w)-k}A_{w_1}(n)\cdots A_{w_k}(n)\\
-\sum_{u\prec w}A_u(n) ,
\end{multline*}
which shows that $A_w(n)+(-1)^{\ell(w)}A_{R(w)}(n)$ can 
always be written in terms of sums of length less than
than $\ell(w)$.  Cf. the discussion in \cite[\S6]{W03}.

\section{EXAMPLE:  SYMMETRIC SUMS}

To illustrate the use of the techniques introduced above, we 
show how symmetric linear combinations of the 
$A_w$ and $S_w$ can be written in terms of ordinary (length 1)
harmonic sums.  
We start in $\E_r$.  Note that the symmetric group
$\Si_k$ acts on words of $k$ letters by permutation, i.e.,
$$
\si\cdot a_1\cdots a_k=a_{\si^{-1}(1)}\cdots a_{\si^{-1}(k)} .
$$
Fix a word $w=a_1a_2\cdots a_k$ of $\E_r$, and for a set partition
$\C=\{C_1,\dots,C_p\}$ of $\{1,2,\dots,k\}$ let
$$
\C(w)=[a_i,i\in C_1]*[a_i,i\in C_2]*\cdots*[a_i,i\in\C_p].
$$
Then repeated use of the multiplication rule (\ref{mult2})
allows $\C(w)$ to be
written as
$$
\sum_{\{B_1,\dots,B_q\}\preceq\C}\sum_{\si\in\Si_q}
\si\cdot [a_i,i\in B_1]\cdots [a_i,i\in B_q],
$$
where $\preceq$ is the partial order given by refinement.  Now
M\"obius inversion can be applied to give
\begin{multline}
\sum_{\si\in\Si_p}\si\cdot[a_i,i\in C_1]\cdots[a_i,i\in C_p]=\\
\sum_{\B=\{B_1,\dots,B_q\}\preceq\C}\mu(\B,\C)\B(w)
\label{sym2prod}
\end{multline}
where $\mu$ is the M\"obius function for the partially ordered set
of partitions of $\{1,2,\dots,k\}$.  When $\C=\{\{1\},\{2\},\dots,\{k\}\}$,
then $\mu(\B,\C)=c(\B)$, where 
$$
c(\B)=(-1)^{k-q}(\card B_1-1)!\cdots(\card B_q-1)!
$$
(see Example 3.10.4 of \cite{S97}). In this case equation 
(\ref{sym2prod}) is
\begin{multline}
\sum_{\si\in\Si_k}\si\cdot w=\\
\sum_{\B=\{B_1,\dots,B_q\}}c(\B)[a_i,i\in B_1]*\cdots*
[a_i,i\in B_q] ,
\label{w2prod}
\end{multline}
where the sum on the right-hand side is over all partitions $\B$ of
$\{1,\dots,k\}$.  Now apply $RS$ to both sides of equation (\ref{w2prod}) 
(and cancel signs) to get
\begin{multline}
\sum_{\si\in\Si_k}\si\cdot\overline w=\\
\sum_{\B=\{B_1,\dots,B_q\}}|c(\B)|[a_i,i\in B_1]*\cdots*
[a_i,i\in B_q] .
\label{barw2prod}
\end{multline}
Finally, we can apply the homomorphism $\rho_n$ to equations (\ref{w2prod}) 
and (\ref{barw2prod}) to get formulas for symmetric combinations of multiple 
harmonic sums in terms of ordinary harmonic sums:
\begin{multline}
\sum_{\si\in\Si_k}A_{\si\cdot w}(n)=\\
\sum_{\B=\{B_1,\dots,B_q\}}c(\B)A_{[a_i,i\in B_1]}(n)\cdots
A_{[a_i,i\in B_q]}(n) 
\label{A2prod}
\end{multline}
\begin{multline}
\sum_{\si\in\Si_k}S_{\si\cdot w}(n)=\\
\sum_{\B=\{B_1,\dots,B_q\}}|c(\B)|A_{[a_i,i\in B_1]}(n)\cdots
A_{[a_i,i\in B_q]}(n) .
\label{S2prod}
\end{multline}
Equations (\ref{A2prod}) and (\ref{S2prod}) generalize Theorem 4.1 of 
\cite{H04} (which
is the case $r=1$).  They may be compared to the corresponding formulas 
for multiple zeta values, which appear as Theorems 2.2 and 2.1, 
respectively, of \cite{H92}.  Equation (\ref{S2prod}) should also
be compared to equations (2.37-2.41) of \cite{B03}, which exhibit
the cases $k=2,\dots,6$ for $r=2$.
\par
In the special case $w=a^k$ (i.e., $w$ is a power of a single letter),
equation (\ref{S2prod}) reduces to
\begin{multline}
S_{a^k}(n)=\\
\frac1{k!}\sum_{\B=\{B_1,\dots,B_q\}}|c(\B)|A_{[b_1a]}(n)\cdots 
A_{[b_qa]}(n) ,
\label{pwr2prod}
\end{multline}
where $b_i=\card B_i$ and $[ka]$ means $[a,a,\dots,a]$ with $k$
repetitions of $a$.
Now for a given (unordered) sequence of 
positive integers $b_1,\dots,b_q$ adding up to $k$, there are 
$$
\frac1{m_1!\cdots m_k!}\frac{k!}{b_1!\cdots b_q!}
$$
partitions of the set $\{1,2\dots,k\}$ having the $b_i$ as block
sizes, where $m_s=\card\{b_i| b_i=s\}$.  
Thus, equation (\ref{pwr2prod}) can be written as
\begin{multline}
S_{a^k}(n)=\\
\sum_{b_1+\dots+b_q=k}\frac1{m_1!\cdots m_k!}\frac1{b_1}A_{[b_1a]}(n)
\cdots \frac1{b_q}A_{[b_qa]}(n) ,
\end{multline}
where the sum is over all integer partitions of $k$
(cf. equations (2.42-2.46) of \cite{B03}).
There is an analogous formula for $A_{a^k}(n)$ differing 
from (21) only in the presence of signs.

\end{document}